\documentclass[12pt]{article}
\usepackage{amssymb}
\usepackage{amsmath}
\usepackage{amsfonts}

\input{tcilatex}
\begin{document}

\begin{center}
\textbf{The canonical projection associated to certain possibly infinite
generalized iterated function system as a fixed point}

\bigskip

\textit{Radu MICULESCU} and \textit{Silviu-Aurelian URZICEANU}

\bigskip
\end{center}

\textbf{Abstract. }{\small In this paper, influenced by the ideas from A.
Mihail, The canonical projection between the shift space of an IIFS and its
attractor as a fixed point, Fixed Point Theory Appl., 2015, Paper No. 75, 15
p., we associate to every generalized iterated function system }$\mathcal{F}$%
{\small \ (of order }$m${\small ) an operator }$H_{\mathcal{F}}:\mathcal{C}%
^{m}\rightarrow \mathcal{C}${\small , where }$\mathcal{C}${\small \ stands
for the space of continuous functions from the shift space on the metric
space corresponding to the system. We provide sufficient conditions (on the
constitutive functions of }$\mathcal{F}${\small ) for the operator }$H_{%
\mathcal{F}}${\small \ to be continuous, contraction, }$\varphi ${\small %
-contraction, Meir-Keeler or contractive. We also give sufficient condition
under which }$H_{\mathcal{F}}${\small \ has a unique fixed point }$\pi _{0}$%
{\small . Moreover, we prove that, under these circumstances, the closer of
the imagine of }$\pi _{0}${\small \ is the attractor of }$\mathcal{F}$%
{\small \ and that }$\pi _{0}${\small \ is the canonical projection
associated to }$\mathcal{F}${\small . In this way we give a partial answer
to the open problem raised on the last paragraph of the above mentioned
Mihail's paper.}

\bigskip

\textit{Key words and phrases:}{\small \ }possibly infinite generalized
iterated function system, canonical projection, attractor, fixed point, $%
\varphi $-contraction, Meir-Keeler function

\textit{2010 Mathematics Subject Classification:} Primary 28A80; Secondary
54H20

\bigskip

\textbf{1. Introduction}

\bigskip

As part of the effort to generalize the concept of iterated function system
introduced by J. Hutchinson (see [3]), R. Miculescu and A. Mihail (see [6]
and [8]) proposed the concept of generalized iterated function system. More
precisely, given $m\in \mathbb{N}$ and a metric space $(X,d)$, a generalized
iterated function system (for short a GIFS) of order $m$ is a finite family
of functions $f_{1},...,f_{n}:X^{m}\rightarrow X$ satisfying certain
contractive conditions. They proved that there exists a unique attractor of
a GIFS and studied some of its properties. F. Strobin (see [14]) proved
that, for any $m\geq 2$, there exists a Cantor subset of the plane which is
an attractor of some GIFS of order $m$, but is not an attractor of a GIFS of
order $m-1$. This shows that GIFSs are real generalizations of iterated
function systems. Certain algorithms generating images of attractors of
GIFSs could be found in [4]. Let us list some extensions of the concept of
GIFS: a) D. Dumitru (see [1] and [2]) investigated generalized iterated
function systems consisting of Meir-Keeler functions; b) F. Strobin and J.
Swaczyna (see [15]) extended the concept of GIFS to the more general setting
of $\varphi $-contractions; c) N. Secelean (see [13]) studied countable
iterated function systems consisting of generalized contraction mappings on
the product space $X^{I}$ into $X$, where $I\subseteq \mathbb{N}$; d) E.
Oliveira and F. Strobin (see [11]) defined the notion of generalized
iterated fuzzy function system. Moreover, the Hutchinson measure associated
with a generalized iterated function system was studied in [7] (for GIFS
with probabilities), in [5] (for generalized iterated function systems with
place dependent probabilities) and in [12].

The canonical projection associated to an iterated function system is a
crucial tool in the study of topological properties of the attractor of such
a system. A significant position from the point of view of this paper is
occupied by [10]. More precisely it is proved there that for a possibly
infinite iterated function system, in two cases (namely: a) the constitutive
functions of the system are uniformly Meir-Keeler; b) the metric space
associated to the system is compact and the system consists of a finite
number of contractive functions), the canonical projection between the shift
space of the system and its attractor can be viewed as a fixed point.

The concept of code space for GIFSs was introduced by A. Mihail (see [9])
and reformulated by F. Strobin and J. Swaczyna (see [16]) in order to treat
the problem of connectedness of the attractor of a GIFS.

In this paper, inspired by the ideas from [10], we associate to a
generalized iterated function system $\mathcal{F}$ (of order $m$) an
operator $H_{\mathcal{F}}:\mathcal{C}^{m}\rightarrow \mathcal{C}$, where $%
\mathcal{C}$ stands for the space of continuous functions from the shift
space on the metric space corresponding to the system. In section 3, we
provide sufficient conditions (on the constitutive functions of $\mathcal{F}$%
) for the operator $H_{\mathcal{F}}$ to be continuous, contraction, $\varphi 
$-contraction, Meir-Keeler or contractive. In section 4, we give sufficient
condition under which $H_{\mathcal{F}}$ has a unique fixed point $\pi _{0}$
(see Theorem 4.1). Moreover, we prove that, under these conditions, the
closer of the imagine of $\pi _{0}$ is the attractor of $\mathcal{F}$ (see
Theorem 4.2) and that $\pi _{0}$ is the canonical projection associated to $%
\mathcal{F}$ (see Theorem 4.3). Our results can be considered as a partial
answer to the open problem raised at the end of [10].

\bigskip

\textbf{2. Preliminaries}

\bigskip

\textbf{A. The Hausdorff-Pompeiu metric}

\bigskip

For a metric space $(X,d)$, the function $h:\mathcal{B}(X)\times \mathcal{B}%
(X)\rightarrow \lbrack 0,\infty )$ given by 
\begin{equation*}
h(A,B)=\max \{d(A,B),d(B,A)\}\text{,}
\end{equation*}%
for every $A,B\in \mathcal{B}(X)$, where $d(A,B)=\underset{x\in A}{\sup }%
\underset{y\in B}{\inf }d(x,y)$, is a metric on $\mathcal{B}(X)$ which is
called the Hausdorff-Pompeiu metric. Here, by $\mathcal{B}(X)$ we mean the
set of all non-empty, closed and bounded subsets of $X$. In the sequel by $%
\mathcal{K}(X)$ we mean the set of all non-empty compact subsets of $X$.

\bigskip

\textbf{B. The metric space }$(X^{m},d_{\max })$

\bigskip

For a metric space $(X,d)$ and $m\in \mathbb{N}^{\ast }$, we endow the
Cartesian product $X^{m}$ with the maximum metric $d_{\max }$ defined by 
\begin{equation*}
d_{\max }((x_{1},...,x_{m}),(y_{1},...,y_{m}))=\max
\{d(x_{1},y_{1}),...,d(x_{m},y_{m})\}\text{,}
\end{equation*}%
for all $(x_{1},...,x_{m}),(y_{1},...,y_{m})\in X^{m}$.

For a metric space $(X,d)$ and $m\in \mathbb{N}^{\ast }$, we define
inductively the spaces $X_{1}$, $X_{2}$, ...., $X_{k}$, ... in the following
way:%
\begin{equation*}
X_{1}=\underset{m\text{ times}}{X\times X\times ...\times X}=X^{m}
\end{equation*}%
and 
\begin{equation*}
X_{k+1}=\underset{m\text{ times}}{X_{k}\times X_{k}\times ...\times X_{k}}
\end{equation*}%
for every $k\in \mathbb{N}^{\ast }$.

We endow $X_{k}$ with the maximum metric for every $k\in \mathbb{N}^{\ast }$%
. Let us to lay stress upon the fact that $X_{k}$ is isometric to $X^{m^{k}}$
with the maximum metric for every $k\in \mathbb{N}^{\ast }$.

\bigskip

\textbf{C.} \textbf{The Mihail-Strobin\&Swaczyna generalized code space}

\bigskip

The notion of code space associated to a generalized iterated function
system was introduced by A. Mihail (see [9]). A different but equivalent
concept which can be easier handled is due to F. Strobin and J. Swaczyna
(see [16]).

\bigskip

Given $m\in \mathbb{N}^{\ast }$ and a set $I$, we define inductively the
sets $\Omega _{1}$, $\Omega _{2}$, ...., $\Omega _{k}$, ... in the following
way:

\begin{equation*}
\Omega _{1}=I
\end{equation*}%
and 
\begin{equation*}
\Omega _{k+1}=\underset{m\text{ times}}{\Omega _{k}\times \Omega _{k}\times
...\times \Omega _{k}}
\end{equation*}
for every $k\in \mathbb{N}^{\ast }$.

We also consider the sets%
\begin{equation*}
\Omega =\Omega _{1}\times \Omega _{2}\times ...\times \Omega _{k}\times ...
\end{equation*}%
and%
\begin{equation*}
_{k}\Omega =\Omega _{1}\times \Omega _{2}\times ...\times \Omega _{k}\text{,}
\end{equation*}%
where $k\in \mathbb{N}^{\ast }$.

For $i\in \{1,2,...,m\}$, $k\in \mathbb{N}$, $k\geq 2$ and $\alpha =\alpha
^{1}\alpha ^{2}...\alpha ^{k}\in $ $_{k}\Omega $, where $\alpha ^{2}=\alpha
_{1}^{2}\alpha _{2}^{2}...\alpha _{m}^{2}\in \Omega _{2}$, ..., $\alpha
^{k}=\alpha _{1}^{k}\alpha _{2}^{k}...\alpha _{m}^{k}\in \Omega _{k}$, we
consider 
\begin{equation*}
\alpha (i)=\alpha _{i}^{2}\alpha _{i}^{3}...\alpha _{i}^{k}\in _{k-1}\Omega 
\text{.}
\end{equation*}

For $\alpha \in \Omega $ and $i\in \{1,2,...,m\}$, we define $\alpha (i)$ in
a similar manner.

\bigskip

\textbf{Definition 1.1.} $\Omega $ \textit{is called the
Mihail-Strobin\&Swaczyna generalized code space.}

\bigskip

\textbf{Remark 1.2.} \textit{Endowed with the metric}\textbf{\ }$d$\textit{\
given by} 
\begin{equation*}
d(\alpha ,\beta )=\underset{k\in \mathbb{N}}{\sum }C^{k}d(\alpha ^{k},\beta
^{k})\text{,}
\end{equation*}%
\textit{for every }$\alpha =\alpha ^{1}\alpha ^{2}...\alpha ^{i}\alpha
^{i+1}...$\textit{, }$\beta =\beta ^{1}\beta ^{2}...\beta ^{i}\beta
^{i+1}...\in \Omega $\textit{, where }$d(\alpha ^{k},\beta ^{k})=\{%
\begin{array}{cc}
1\text{,} & \alpha ^{k}\neq \beta ^{k} \\ 
0\text{,} & \alpha ^{k}=\beta ^{k}%
\end{array}%
$\textit{\ and} $C\in (0,1)$, $(\Omega ,d)$\textit{\ becomes a complete
metric space.}

\bigskip

\textbf{Remark 1.3.} \textit{If }$I$\textit{\ is finite, then the metric
space }$(\Omega ,d)$\textit{\ is compact.}

\bigskip

\textbf{D. Generalized possibly infinite iterated function systems}

\bigskip

\textbf{Definition 1.4.} \textit{A generalized possibly infinite iterated
function system of order }$m\in \mathbb{N}^{\ast }$\textit{\ is a pair }$%
\mathcal{F}=((X,d),(f_{i})_{i\in I})$\textit{, where }$(X,d)$ \textit{is a
metric space, }$f_{i}:X^{m}\rightarrow X$ \textit{is continuous for every} $%
i\in I$\textit{\ and the family of functions }$(f_{i})_{i\in I}$\textit{\ is
bounded} \textit{(i.e. }$\underset{_{i\in I}}{\cup }f_{i}(B)$ \textit{is
bounded for each bounded subset }$B$ \textit{of }$X^{m}$\textit{).}

\textit{The function }$\mathcal{F}_{\mathcal{F}}:(\mathcal{B}%
(X))^{m}\rightarrow \mathcal{B}(X)$, \textit{given by} 
\begin{equation*}
\mathcal{F}_{\mathcal{F}}(B_{1},...,B_{m})=\overline{\underset{i\in I}{\cup }%
f_{i}(B_{1},...,B_{m})}\text{,}
\end{equation*}%
\textit{for all} $(B_{1},...,B_{m})\in (\mathcal{B}(X))^{m}$, \textit{is
called the fractal operator associated to} $\mathcal{F}$.

\textit{If }$\mathcal{F}_{\mathcal{F}}$ \textit{has a unique fixed point,
then it is called the attractor of }$\mathcal{F}$ \textit{and it is denoted
by} $A_{\mathcal{F}}$.

\bigskip

If the set $I$ is finite, then $\mathcal{F}_{\mathcal{F}}((\mathcal{K}%
(X))^{m})\subseteq \mathcal{K}(X)$ and we make the convention to still
denote the function $(B_{1},...,B_{m})\in (\mathcal{K}(X))^{m}\rightarrow 
\mathcal{F}_{\mathcal{F}}(B_{1},...,B_{m})\in \mathcal{K}(X)$ by $\mathcal{F}%
_{\mathcal{F}}$. In this case $A_{\mathcal{F}}\in \mathcal{K}(X)$.

\bigskip

For a generalized possibly infinite iterated function system $\mathcal{F}%
=((X,d),(f_{i})_{i\in I})$ of order $m$, we define inductively a family of
functions $\{f_{\alpha }:X_{k}\rightarrow X\mid \alpha \in $ $_{k}\Omega \}$
for every $k\in \mathbb{N}^{\ast }$ in the following way:

i) For $k=1$, the family is $(f_{i})_{i\in I}$.

ii) If the functions $f_{\alpha }$, where $\alpha \in $ $_{k}\Omega $, have
been defined, then, we define 
\begin{equation*}
f_{\alpha }(x_{1},x_{2},...,x_{m})=f_{\alpha ^{1}}(f_{\alpha
(1)}(x_{1}),...,f_{\alpha (m)}(x_{m}))
\end{equation*}%
for every $\alpha =\alpha ^{1}\alpha ^{2}...\alpha ^{k}\alpha ^{k+1}\in $ $%
_{k+1}\Omega $, where $\alpha ^{1}\in \Omega _{1}$, $\alpha ^{2}\in \Omega
_{2}$, $...$, $\alpha ^{k}\in \Omega _{k}$, $\alpha ^{k+1}\in \Omega _{k+1}$%
, $(x_{1},x_{2},...,x_{m})\in X_{k+1}=\underset{m\text{ times}}{X_{k}\times
X_{k}\times ...\times X_{k}}$.

Note that the above introduced families of functions are natural
generalizations of compositions of functions since if $m=1$, then $%
_{k}\Omega =I^{k}$ and if $\omega =\omega ^{1}\omega ^{2}...\omega ^{k}\in $ 
$_{k}\Omega $, then $f_{\omega }=f_{\omega ^{1}}\circ ...\circ f_{\omega
^{k}}$.

\bigskip

\textbf{E. The operator }$H_{\mathcal{F}}$ \textbf{associated to a
generalized possibly infinite iterated function system}

\bigskip

For a generalized possibly infinite iterated function system $\mathcal{F}%
=((X,d),(f_{i})_{i\in I})$ of order $m$, we consider the operator $H_{%
\mathcal{F}}:\mathcal{C}^{m}\rightarrow \mathcal{C}$ given by%
\begin{equation*}
H_{\mathcal{F}}(g_{1},...,g_{m})(\alpha )=f_{\alpha ^{1}}(g_{1}(\alpha
(1)),...,g_{m}(\alpha (m)))\text{,}
\end{equation*}%
for every $g_{1},...,g_{m}\in \mathcal{C}$ and every $\alpha =\alpha
^{1}\alpha ^{2}...\alpha ^{k}...\in \Omega $, $\alpha ^{k}\in \Omega _{k}$
for every $k\in \mathbb{N}^{\ast }$, where the metric space $(\mathcal{C}%
,d_{u})$ is described by 
\begin{equation*}
\mathcal{C=\{}f:\Omega \rightarrow X\mid f\text{ is continuous and bounded}\}
\end{equation*}%
and 
\begin{equation*}
d_{u}(f,g)=\underset{\alpha \in \Omega }{\sup }\text{ }d(f(\alpha ),g(\alpha
))
\end{equation*}%
for every $f,g\in \mathcal{C}$.

\bigskip

\textbf{Remark 1.5}.

i) \textit{If} $(X,d)$ \textit{is complete, then} $\mathcal{C}$ \textit{is
complete.}

ii) $H_{\mathcal{F}}$ \textit{is well defined}, i.e. $H_{\mathcal{F}%
}(g_{1},...,g_{m})$ is continuous and bounded for all $g_{1},...,g_{m}\in 
\mathcal{C}$. Indeed, on one hand, the continuity follows from the following
facts: $\Omega =\underset{i\in I}{\cup }\Omega ^{i}$, where $\Omega
^{i}=\{\alpha =\alpha ^{1}\alpha ^{2}...\alpha ^{i}\alpha ^{i+1}...\in
\Omega \mid \alpha ^{1}=i\}$, and the restriction of $H_{\mathcal{F}%
}(g_{1},...,g_{m})$ to the open set $\Omega ^{i}$ is continuous for every $%
i\in I$. On the other hand, the boundedness follows from the boundedness of
the family of functions $(f_{i})_{i\in I}$, the boundedness of the functions 
$g_{1},...,g_{m}$ and from the fact that%
\begin{equation*}
H_{\mathcal{F}}(g_{1},...,g_{m})(\Omega )=H_{\mathcal{F}}(g_{1},...,g_{m})(%
\underset{i\in I}{\cup }\Omega ^{i})=
\end{equation*}%
\begin{equation*}
=\underset{i\in I}{\cup }H_{\mathcal{F}}(g_{1},...,g_{m})(\Omega ^{i})=%
\underset{i\in I}{\cup }f_{i}(g_{1}(\Omega ),...,g_{m}(\Omega ))\text{.}
\end{equation*}

\bigskip

\textbf{F.} \textbf{Some classes of functions }$f:X^{m}\rightarrow X$\textbf{%
\ and their fixed points}

\bigskip

Given a set $X$,\textit{\ }$m\in \mathbb{N}^{\ast }$ and a function $%
f:X^{m}\rightarrow X$, we define inductively a family of functions $%
f^{[k]}:X^{m^{k}}\rightarrow X$, $k\in \mathbb{N}^{\ast }$, in the following
way:

i)%
\begin{equation*}
f^{[1]}=f
\end{equation*}

ii) assuming that we have defined $f^{[k]}$, then 
\begin{equation*}
f^{[k+1]}(x_{1},...,x_{m})=f(f^{[k]}(x_{1}),...,f^{[k]}(x_{m}))\text{,}
\end{equation*}%
for every $(x_{1},...,x_{m})\in \underset{m\text{ times}}{X^{m^{k}}\times
...\times X^{m^{k}}}=X^{m^{k+1}}=X_{k+1}$.

Note that for $m=1$, we have $f^{[k]}=\underset{k\text{ times}}{f\circ
...\circ f}$.

\bigskip

\textbf{Definition 1.6.} \textit{Given a set} $X$ \textit{and }$m\in \mathbb{%
N}^{\ast }$\textit{, an element }$x$\textit{\ of }$X$\textit{\ is called a
fixed point of a function }$f:X^{m}\rightarrow X$\textit{\ if }%
\begin{equation*}
f(x,...,x)=x\text{\textit{.}}
\end{equation*}

\bigskip

\textbf{Definition 1.7.} \textit{A function }$\varphi :[0,\infty
)\rightarrow \lbrack 0,\infty )$\textit{\ is called a comparison function if
it satisfies the following properties:}

\textit{i) it is nondecreasing;}

\textit{ii) it is right-continuous;}

\textit{iii) }$\varphi (t)<t$\textit{\ for every }$t>0$\textit{.}

\bigskip

\textbf{Definition 1.8. }\textit{Given a metric space }$(X,d)$, $m\in 
\mathbb{N}^{\ast }$\textit{\ and a comparison function }$\varphi $\textit{,
a function }$f:X^{m}\rightarrow X$\textit{\ is called a }$\varphi $\textit{%
-contraction if}%
\begin{equation*}
d(f(x),f(y))\leq \varphi (d_{\max }(x,y))\text{,}
\end{equation*}%
\textit{for all }$x,y\in X^{m}$\textit{.}

\bigskip

\textbf{Definition 1.9.} \textit{Given a metric space} $(X,d)$ \textit{and }$%
m\in \mathbb{N}^{\ast }$\textit{, a function }$f:X^{m}\rightarrow X$\textit{%
\ is called Meir-Keeler if for every }$\varepsilon >0$\textit{\ there exists 
}$\delta _{\varepsilon }>0$\textit{\ such that}%
\begin{equation*}
d(f(x),f(y))<\varepsilon \text{,}
\end{equation*}%
\textit{for all }$x,y\in X^{m}$\textit{\ having the property that }$d_{\max
}(x,y)<\varepsilon +\delta _{\varepsilon }$.

\bigskip

\textbf{Definition 1.10. }\textit{Given a metric space }$(X,d)$\ \textit{and}
$m\in \mathbb{N}^{\ast }$\textit{, a function }$f:X^{m}\rightarrow X$\textit{%
\ is called contractive if}%
\begin{equation*}
d(f(x),f(y))<d_{\max }(x,y)\text{,}
\end{equation*}%
\textit{for all }$x,y\in X^{m}$\textit{, }$x\neq y$.

\bigskip

\textbf{Theorem 1.11 }(see Theorem 19 from\ [2])\textbf{.} \textit{Given a
complete metric space} $(X,d)\ $\textit{and} $m\in \mathbb{N}^{\ast }$%
\textit{, for each Meir-Keeler function }$f:X^{m}\rightarrow X$ \textit{%
there exists a unique fixed point }$x_{0}$\textit{\ of }$f$\textit{\ and the
sequence }$(f^{[n]}(x,...,x))_{n\in \mathbb{N}^{\ast }}$\textit{\ converges
to }$x_{0}$\textit{\ for every }$x\in X$\textit{.}

\bigskip

The following definition is inspired by Definition 2.3 from [10] and
Definition 17 from [2].

\bigskip

\textbf{Definition 1.12.} \textit{Given a metric space} $(X,d)$ \textit{and }%
$m\in \mathbb{N}^{\ast }$\textit{, a family of functions }$%
f_{i}:X^{m}\rightarrow X$\textit{, }$i\in I$, \textit{is called uniformly
Meir-Keeler if for every }$\varepsilon >0$\textit{\ there exist }$\delta
_{\varepsilon },\lambda _{\varepsilon }>0$\textit{\ such that}%
\begin{equation*}
d(f_{i}(x),f_{i}(y))<\varepsilon -\lambda _{\varepsilon }\text{,}
\end{equation*}%
\textit{for all }$i\in I$ \textit{and all }$x,y\in X^{m}$\textit{\ having
the property that }$d_{\max }(x,y)<\varepsilon +\delta _{\varepsilon }$%
\textit{.}

\bigskip

\textbf{G. Special classes of generalized iterated function systems}

\bigskip

\textbf{Theorem 1.13 }(see Theorem 3.11 from [15])\textbf{.} \textit{For
each generalized possibly infinite iterated function system }$\mathcal{F}%
=((X,d),(f_{i})_{i\in I})$\textit{, where }$(X,d)$ \textit{is a complete
metric space, }$I$ \textit{is finite and all the functions }$f_{i}$\textit{\
are }$\varphi $\textit{-contractions for some comparison function }$\varphi $%
\textit{, there exists a unique set }$A_{\mathcal{F}}\in \mathcal{K}(X)$ 
\textit{such that} $\mathcal{F}_{\mathcal{F}}(A_{\mathcal{F}},...,A_{%
\mathcal{F}})=A_{\mathcal{F}}$\textit{\ }(i.e. $\mathcal{F}$ has attractor%
\textit{).}

\bigskip

\textbf{Remark 1.14}.

i) In the framework of the above theorem, \textit{the set }$\underset{k\in 
\mathbb{N}}{\cap }f_{\alpha 1...\alpha ^{k}}(A_{\mathcal{F}})$\textit{\
consists on a single element} denoted by $x_{\alpha }$ for every $\alpha
=\alpha ^{1}...\alpha ^{i}...\in \Omega $. The function $\pi :\Omega
\rightarrow X$ given by $\pi (\alpha )=x_{\alpha }$, for every $\alpha \in
\Omega $, is called \textit{the canonical projection associated to} $%
\mathcal{F}$. For the properties of this function see Theorem 3.7, Corollary
3.9 and Theorem 3.11 from [16].

ii) The same line of reasoning used in [15] and [16] leads to the following
conclusion: \textit{Each} \textit{generalized possibly infinite iterated
function system }$\mathcal{F}=((X,d),(f_{i})_{i\in I})$\textit{, where }$%
(X,d)$ \textit{is a complete metric space and all the functions }$f_{i}$%
\textit{\ are }$\varphi $\textit{-contractions for some comparison function }%
$\varphi $\textit{, has attractor}, i.e. there exists a unique set $A_{%
\mathcal{F}}\in \mathcal{B}(X)$ such that $\mathcal{F}_{\mathcal{F}}(A_{%
\mathcal{F}},...,A_{\mathcal{F}})=A_{\mathcal{F}}$\textit{. }The function%
\textit{\ }$\pi :\Omega \rightarrow X$, described by $\{\pi (\alpha )\}=%
\underset{k\in \mathbb{N}}{\cap }\overline{f_{\alpha ^{1}...\alpha ^{k}}(A_{%
\mathcal{F}})}$ for every $\alpha \in \Omega $, which is called \textit{the
canonical projection associated to} $\mathcal{F}$, has the property that $A_{%
\mathcal{F}}=\overline{\pi (\Omega )}$.

\bigskip

\textbf{Theorem 1.15 }(see Theorem 32 from [2])\textbf{.} \textit{For each
generalized possibly infinite iterated function system }$\mathcal{F}%
=((X,d),(f_{i})_{i\in I})$\textit{, where }$(X,d)$ \textit{is a complete
metric space and the family of functions }$(f_{i})_{i\in I}$ \textit{is
uniformly Meir-Keeler, there exists a unique set }$A_{\mathcal{F}}\in 
\mathcal{B}(X)$ \textit{such that} $\mathcal{F}_{\mathcal{F}}(A_{\mathcal{F}%
},...,A_{\mathcal{F}})=A_{\mathcal{F}}$\textit{\ }(i.e.\textit{\ }$\mathcal{F%
}$ has attractor\textit{).}

\bigskip

\textbf{Remark 1.16}. \textit{The considerations from Remark 1.14 are also
valid in the framework of Theorem 1.15. }The arguments supporting this
Remark are almost the same with the ones used for Remark 1.14. The only fact
which needs a special attention is the justification of the following
equality:

\begin{equation}
\underset{n\rightarrow \infty }{\lim }\underset{\alpha \in \text{ }%
_{n}\Omega }{\sup }\text{diam }(f_{\alpha }(A_{\mathcal{F}},...,A_{\mathcal{F%
}}))=0\text{,}  \tag{1}
\end{equation}%
where by diam $(A)$ we mean the diameter of the subset $A$ of $X$. In order
to prove $(1)$, we adopt the following notation: $\underset{\alpha \in \text{
}_{n}\Omega }{\sup }$diam $(f_{\alpha }(A_{\mathcal{F}},...,A_{\mathcal{F}}))%
\overset{not}{=}d_{n}$, $n\in \mathbb{N}^{\ast }$. Let us note that $%
(d_{n})_{n\in \mathbb{N}^{\ast }}$ is a decreasing sequence of positive real
numbers (since $f_{\alpha ^{1}\alpha ^{2}...\alpha ^{n}\alpha ^{n+1}}(A_{%
\mathcal{F}},...,A_{\mathcal{F}})\subseteq f_{\alpha ^{1}\alpha
^{2}...\alpha ^{n}}(A_{\mathcal{F}},...,A_{\mathcal{F}})$ for every $\alpha
^{1}\in \Omega _{1}$, $\alpha ^{2}\in \Omega _{2}$, ..., $\alpha ^{n}\in
\Omega _{n}$, $\alpha ^{n+1}\in \Omega _{n+1}$, $n\in \mathbb{N}^{\ast }$),
so there exists $l\geq 0$ such that $\underset{n\rightarrow \infty }{\lim }%
d_{n}=l$. The justification of $(1)$ is done if we prove that $l=0$. Let us
suppose, by reductio ad absurdum, that $l>0$. Then, in view of the fact that
the family of functions $(f_{i})_{i\in I}$ is uniformly Meir-Keeler, there
exist\textit{\ }$\delta _{l},\lambda _{l}>0$\textit{\ }such that 
\begin{equation}
d(f_{i}(x),f_{i}(y))<l-\lambda _{l}\text{,}  \tag{2}
\end{equation}%
for all $i\in I$ and all $x,y\in X^{m}$\ having the property that $d_{\max
}(x,y)<l+\delta _{l}$. Now, based on the fact that $\underset{n\rightarrow
\infty }{\lim }d_{n}=l$, let us choose $n_{0}\in \mathbb{N}^{\ast }$ such
that $d_{n_{0}}<l+\delta _{l}$ and note that 
\begin{equation*}
d(f_{\alpha }(x_{1},...,x_{m}),f_{\alpha }(y_{1},...,y_{m}))=
\end{equation*}%
\begin{equation}
=d(f_{\alpha ^{1}}((f_{\alpha (1)}(x_{1}),...,f_{\alpha
(m)}(x_{m}))),f_{\alpha ^{1}}((f_{\alpha (1)}(y_{1}),...,f_{\alpha
(m)}(y_{m}))))\overset{(2)}{<}l-\lambda _{l}\text{,}  \tag{3}
\end{equation}%
since%
\begin{equation*}
d_{\max }((f_{\alpha (1)}(x_{1}),...,f_{\alpha (m)}(x_{m})),(f_{\alpha
(1)}(y_{1}),...,f_{\alpha (m)}(y_{m})))=
\end{equation*}%
\begin{equation*}
=\max \{d(f_{\alpha (1)}(x_{1}),f_{\alpha (1)}(y_{1})),...,d(f_{\alpha
(m)}(x_{m}),f_{\alpha (m)}(y_{m}))\}\leq
\end{equation*}%
\begin{equation*}
\leq \max \{\text{diam }f_{_{\alpha (1)}}(A_{\mathcal{F}},...,A_{\mathcal{F}%
}),...,\text{diam }f_{_{\alpha (m)}}(A_{\mathcal{F}},...,A_{\mathcal{F}%
})\}\leq d_{n_{0}}<l+\delta _{l}\text{,}
\end{equation*}%
for every $\alpha =\alpha ^{1}\alpha ^{2}...\alpha ^{n_{0}}\alpha
^{n_{0}+1}\in $ $_{n_{0}+1}\Omega $, where $\alpha ^{1}\in \Omega _{1}$, $%
\alpha ^{2}\in \Omega _{2}$, ..., \linebreak $\alpha ^{n_{0}}\in \Omega
_{n_{0}}$, $\alpha ^{n_{0}+1}\in \Omega _{n_{0}+1}$, and every $%
x_{1},...,x_{m},y_{1},...,y_{m}\in (A_{\mathcal{F}})^{m^{n_{0}}}$.
Consequently, we obtain the following contradiction: $l\leq d_{n_{0}+1}%
\overset{(3)}{\leq }l-\lambda _{l}<l$.

\bigskip

\textbf{3. The properties of the operator }$H_{\mathcal{F}}$

\bigskip

In this section we present some results which give sufficient conditions (on
the constitutive functions of a generalized possibly infinite iterated
function system $\mathcal{F}$) for the operator $H_{\mathcal{F}}$ to be
continuous, generalized contraction, generalized $\varphi $-contraction,
Meir-Keeler or contractive.

\bigskip

\textbf{Proposition 3.1.} \textit{For every generalized possibly infinite
iterated function system} $\mathcal{F}=((X,d),(f_{i})_{i\in I})$ \textit{of
order }$m$\textit{\ such that the family of functions }$(f_{i})_{i\in I}$%
\textit{\ is uniformly equicontinuous (}i.e. for each $\varepsilon >0$ there
exists $\delta _{\varepsilon }>0$ such that for every $i\in I$ and every $%
x,y\in X^{m}$ having the property that $d_{\max }(x,y)<\delta _{\varepsilon
} $ we have $d(f_{i}(x),f_{i}(y))<\varepsilon $\textit{), the operator }$H_{%
\mathcal{F}}$ \textit{is continuous.}

\textit{Proof}. We are going to prove that for each sequence $(g_{n})_{n\in 
\mathbb{N}}$ of elements from $\mathcal{C}^{m}$ and $g\in \mathcal{C}^{m}$
such that $\underset{n\rightarrow \infty }{\lim }g_{n}=g$, we have $\underset%
{n\rightarrow \infty }{\lim }H_{\mathcal{F}}(g_{n})=H_{\mathcal{F}}(g)$. Let
us suppose that $g_{n}=(g_{n}^{1},...,g_{n}^{m})$ and $g=(g^{1},...,g^{m})$,
where $g_{n}^{1},...,g_{n}^{m},g^{1},...,g^{m}\in \mathcal{C}$.

Let us fix $\varepsilon >0$.

Since the family of functions $(f_{i})_{i\in I}$ is uniformly
equicontinuous, there exists $\delta _{\varepsilon }>0$ such that%
\begin{equation}
d(f_{i}(x),f_{i}(y))<\varepsilon \text{,}  \tag{1}
\end{equation}%
for every $i\in I$ and every $x,y\in X^{m}$ having the property that $%
d_{\max }(x,y)<\delta _{\varepsilon }$.

Since $\underset{n\rightarrow \infty }{\lim }g_{n}=g$, there exists $%
n_{\varepsilon }\in \mathbb{N}$ such that 
\begin{equation*}
d_{\max }(g_{n},g)=\max
\{d_{u}(g_{n}^{1},g^{1}),...,d_{u}(g_{n}^{m},g^{m})\}<\delta _{\varepsilon }%
\text{,}
\end{equation*}%
for every $n\in \mathbb{N}$, $n\geq n_{\varepsilon }$. Consequently 
\begin{equation*}
d(g_{n}^{i}(\alpha ),g^{i}(\alpha ))\leq \underset{\alpha \in \Omega }{\sup }%
\text{ }d(g_{n}^{i}(\alpha ),g^{i}(\alpha ))=d_{u}(g_{n}^{i},g^{i})\leq
d_{\max }(g_{n},g)<\delta _{\varepsilon }\text{,}
\end{equation*}%
for all $i\in \{1,...,m\}$, $\alpha \in \Omega $ and $n\in \mathbb{N}$, $%
n\geq n_{\varepsilon }$. Therefore%
\begin{equation}
d_{\max }((g_{n}^{1}(\alpha (1)),...,g_{n}^{m}(\alpha (m)),(g^{1}(\alpha
(1)),...,g^{m}(\alpha (m))))<\delta _{\varepsilon }\text{,}  \tag{2}
\end{equation}%
for all $\alpha \in \Omega $ and $n\in \mathbb{N}$, $n\geq n_{\varepsilon }$.

Based on $(1)$ and $(2)$ we deduce that 
\begin{equation*}
d(f_{\alpha ^{1}}(g_{n}^{1}(\alpha (1)),...,g_{n}^{m}(\alpha (m))),f_{\alpha
^{1}}(g^{1}(\alpha (1)),...,g^{m}(\alpha (m))))<\varepsilon \text{,}
\end{equation*}%
for all $n\in \mathbb{N}$, $n\geq n_{\varepsilon }$ and all $\alpha =\alpha
^{1}\alpha ^{2}...\alpha ^{k}\alpha ^{k+1}...\in \Omega $, where $\alpha
^{k}\in \Omega _{k}$ for every $k\in \mathbb{N}$. The last inequality takes
the following shape: 
\begin{equation*}
d(H_{\mathcal{F}}(g_{n}^{1},...,g_{n}^{m})(\alpha ),H_{\mathcal{F}%
}(g^{1},...,g^{m})(\alpha ))<\varepsilon \text{,}
\end{equation*}%
for all $\alpha \in \Omega $ and $n\in \mathbb{N}$, $n\geq n_{\varepsilon }$%
. Hence%
\begin{equation*}
d_{u}(H_{\mathcal{F}}(g_{n}^{1},...,g_{n}^{m}),H_{\mathcal{F}%
}(g^{1},...,g^{m}))=
\end{equation*}%
\begin{equation*}
=\underset{\alpha \in \Omega }{\sup }\text{ }d(H_{\mathcal{F}%
}(g_{n}^{1},...,g_{n}^{m})(\alpha ),H_{\mathcal{F}}(g^{1},...,g^{m})(\alpha
))\leq \varepsilon \text{,}
\end{equation*}%
for every $n\in \mathbb{N}$, $n\geq n_{\varepsilon }$, i.e. $\underset{%
n\rightarrow \infty }{\lim }H_{\mathcal{F}}(g_{n})=H_{\mathcal{F}}(g)$. $%
\square $

\bigskip

\textbf{Proposition 3.2.} \textit{For every generalized possibly infinite
iterated function system} $\mathcal{F}=((X,d),(f_{i})_{i\in I})$ \textit{of
order }$m$\textit{\ we have }lip$(H_{\mathcal{F}})\leq \underset{i\in I}{%
\sup }$ lip$(f_{i})$. \textit{In particular, if} $\underset{i\in I}{\sup }$
lip$(f_{i})<1$\textit{, then the operator }$H_{\mathcal{F}}$ \textit{is a
contraction.}

\textit{Proof}. With the notation $\underset{i\in I}{\sup }$ lip$(f_{i})%
\overset{not}{=}C$, we have%
\begin{equation*}
d_{u}(H_{\mathcal{F}}(g),H_{\mathcal{F}}(h))=\text{ }\underset{\alpha \in
\Omega }{\sup }\text{ }d(H_{\mathcal{F}}(g^{1},...,g^{m})(\alpha ),H_{%
\mathcal{F}}(h^{1},...,h^{m})(\alpha ))=
\end{equation*}%
\begin{equation*}
=\underset{\alpha =\alpha ^{1}\alpha ^{2}...\alpha ^{k}\alpha ^{k+1}..\in
\Omega }{\sup }\text{ }d(f_{\alpha ^{1}}(g^{1}(\alpha (1)),...,g^{m}(\alpha
(m))),f_{\alpha ^{1}}(h^{1}(\alpha (1)),...,h^{m}(\alpha (m))))\leq
\end{equation*}%
\begin{equation*}
\leq C\text{ }\underset{\alpha \in \Omega }{\sup }\text{ }d_{\max
}((g^{1}(\alpha (1)),...,g^{m}(\alpha (m))),(h^{1}(\alpha
(1)),...,h^{m}(\alpha (m))))=
\end{equation*}%
\begin{equation*}
=C\text{ }\underset{\alpha \in \Omega }{\sup }\max \{d(g^{1}(\alpha
(1)),h^{1}(\alpha (1))),...,d(g^{m}(\alpha (m)),h^{m}(\alpha (m)))\}\leq
\end{equation*}%
\begin{equation*}
\leq C\text{ }\max \{\underset{\alpha \in \Omega }{\sup }\text{ }%
d(g^{1}(\alpha ),h^{1}(\alpha )),...,\underset{\alpha \in \Omega }{\sup }%
\text{ }d(g^{m}(\alpha ),h^{m}(\alpha ))\}=
\end{equation*}%
\begin{equation*}
=C\text{ }\max \{d_{u}(g^{1},h^{1}),...,d_{u}(g^{m},h^{m})\}=C\text{ }%
d_{\max }(g,h)\text{,}
\end{equation*}%
for every $g=(g^{1},...,g^{m}),h=(h^{1},...,h^{m})\in \mathcal{C}^{m}$. $%
\square $

\bigskip

\textbf{Proposition 3.3.} \textit{For every comparison function }$\varphi $ 
\textit{and every generalized possibly infinite iterated function system} $%
\mathcal{F}=((X,d),(f_{i})_{i\in I})$ \textit{of order }$m$\textit{\ such
that all the functions }$f_{i}$\textit{\ are }$\varphi $\textit{%
-contractions, the operator }$H_{\mathcal{F}}$ \textit{is a }$\varphi $%
\textit{-contraction.}

\textit{Proof}. We have%
\begin{equation*}
d_{u}(H_{\mathcal{F}}(g),H_{\mathcal{F}}(h))=\text{ }\underset{\alpha \in
\Omega }{\sup }\text{ }d(H_{\mathcal{F}}(g^{1},...,g^{m})(\alpha ),H_{%
\mathcal{F}}(h^{1},...,h^{m})(\alpha ))=
\end{equation*}%
\begin{equation*}
=\underset{\alpha =\alpha ^{1}\alpha ^{2}...\alpha ^{k}\alpha ^{k+1}..\in
\Omega }{\sup }\text{ }d(f_{\alpha ^{1}}(g^{1}(\alpha (1)),...,g^{m}(\alpha
(m))),f_{\alpha ^{1}}(h^{1}(\alpha (1)),...,h^{m}(\alpha (m))))\leq
\end{equation*}%
\begin{equation*}
\leq \text{ }\underset{\alpha \in \Omega }{\sup }\text{ }\varphi (d_{\max
}((g^{1}(\alpha (1)),...,g^{m}(\alpha (m))),(h^{1}(\alpha
(1)),...,h^{m}(\alpha (m)))))=
\end{equation*}%
\begin{equation*}
=\text{ }\underset{\alpha \in \Omega }{\sup }\text{ }\varphi (\max
\{d(g^{1}(\alpha (1)),h^{1}(\alpha (1))),...,d(g^{m}(\alpha
(m)),h^{m}(\alpha (m)))\})\leq
\end{equation*}%
\begin{equation*}
\leq \varphi (\max \{\underset{\alpha \in \Omega }{\sup }\text{ }%
d(g^{1}(\alpha ),h^{1}(\alpha )),...,\underset{\alpha \in \Omega }{\sup }%
\text{ }d(g^{m}(\alpha ),h^{m}(\alpha )))\}=
\end{equation*}%
\begin{equation*}
=\varphi (\max \{d_{u}(g^{1},h^{1}),...,d_{u}(g^{m},h^{m})\})=\varphi
(d_{\max }(g,h))\text{,}
\end{equation*}%
for every $g=(g^{1},...,g^{m}),h=(h^{1},...,h^{m})\in \mathcal{C}^{m}$. $%
\square $

\bigskip

\textbf{Proposition 3.4.} \textit{For every generalized possibly infinite
iterated function system} $\mathcal{F}=((X,d),(f_{i})_{i\in I})$ \textit{of
order }$m$\textit{\ such that the family of functions }$(f_{i})_{i\in I}$%
\textit{\ is uniformly Meir-Keeler, the operator }$H_{\mathcal{F}}$ \textit{%
is Meir-Keeler.}

\textit{Proof}. Let $\varepsilon >0$ be fixed, but arbitrarily chosen.

Since the family of functions $(f_{i})_{i\in I}$\ is uniformly Meir-Keeler,
there exist $\delta _{\varepsilon },\lambda _{\varepsilon }>0$\ such that 
\begin{equation}
d(f_{i}(x),f_{i}(y))<\varepsilon -\lambda _{\varepsilon }\text{,}  \tag{1}
\end{equation}%
for all $i\in I$ and all $x,y\in X^{m}$\ having the property that $d_{\max
}(x,y)<\varepsilon +\delta _{\varepsilon }$.

If $g=(g^{1},...,g^{m}),h=(h^{1},...,h^{m})\in \mathcal{C}^{m}$ are such
that $d_{\max }(g,f)<\varepsilon +\delta _{\varepsilon }$, then%
\begin{equation*}
d_{\max }((g^{1}(\alpha (1)),...,g^{m}(\alpha (m))),(h^{1}(\alpha
(1)),...,h^{m}(\alpha (m))))=
\end{equation*}%
\begin{equation*}
=\max \{d(g^{1}(\alpha (1)),h^{1}(\alpha (1))),...,d(g^{m}(\alpha
(m)),h^{m}(\alpha (m)))\}\leq 
\end{equation*}%
\begin{equation*}
\leq \max \{d_{u}(g^{1},h^{1}),...,d_{u}(g^{m},h^{m})\}=d_{\max
}(g,f)<\varepsilon +\delta _{\varepsilon }
\end{equation*}%
for all $\alpha \in \Omega $. Then, taking into account $(1)$, we get 
\begin{equation}
d(f_{\alpha ^{1}}(g^{1}(\alpha (1)),...,g^{m}(\alpha (m))),f_{\alpha
^{1}}(h^{1}(\alpha (1)),...,h^{m}(\alpha (m))))<\varepsilon -\lambda
_{\varepsilon }\text{,}  \tag{2}
\end{equation}%
for all $\alpha =\alpha ^{1}\alpha ^{2}...\alpha ^{k}\alpha ^{k+1}...\in
\Omega $ and all $g=(g^{1},...,g^{m}),h=(h^{1},...,h^{m})\in \mathcal{C}^{m}$
such that $d_{\max }(g,f)<\varepsilon +\delta _{\varepsilon }$. Hence 
\begin{equation*}
\underset{\alpha \in \Omega }{\sup }\text{ }d_{u}(H_{\mathcal{F}%
}(g^{1},...,g^{m})(\alpha ),H_{\mathcal{F}}(h^{1},...,h^{m})(\alpha ))=
\end{equation*}%
\begin{equation*}
=\underset{\alpha \in \Omega }{\sup }\text{ }d(f_{\alpha ^{1}}(g^{1}(\alpha
(1)),...,g^{m}(\alpha (m))),f_{\alpha ^{1}}(h^{1}(\alpha
(1)),...,h^{m}(\alpha (m))))\overset{(2)}{\leq }
\end{equation*}%
\begin{equation*}
\leq \varepsilon -\lambda _{\varepsilon }<\varepsilon \text{,}
\end{equation*}%
for all $g=(g^{1},...,g^{m}),h=(h^{1},...,h^{m})\in \mathcal{C}^{m}$ such
that $d_{\max }(g,f)<\varepsilon +\delta _{\varepsilon }$, where $\alpha
=\alpha ^{1}\alpha ^{2}...\alpha ^{k}\alpha ^{k+1}...$ .

Consequently, as 
\begin{equation*}
d_{u}(H_{\mathcal{F}}(g),H_{\mathcal{F}}(h))=\underset{\alpha \in \Omega }{%
\sup }\text{ }d_{u}(H_{\mathcal{F}}(g^{1},...,g^{m})(\alpha ),H_{\mathcal{F}%
}(h^{1},...,h^{m})(\alpha ))\text{,}
\end{equation*}%
we obtained that $d_{u}(H_{\mathcal{F}}(g),H_{\mathcal{F}}(h))<\varepsilon $
for every $g,h\in \mathcal{C}^{m}$ such that $d_{\max }(g,f)<\varepsilon
+\delta _{\varepsilon }$, i.e. $H_{\mathcal{F}}$ is Meir-Keeler. $\square $

\bigskip

\textbf{Proposition 3.5.} \textit{For every generalized possibly infinite
iterated function system} $\mathcal{F}=((X,d),(f_{i})_{i\in I})$ \textit{of
order }$m$\textit{\ such that }$I$ \textit{is finite and all the functions} $%
f_{i}$\textit{\ are Meir-Keeler, the operator }$H_{\mathcal{F}}$ \textit{is
Meir-Keeler.}

\textit{Proof}. Let $\varepsilon >0$ be a fixed, but arbitrarily chosen.

Since all the functions $f_{i}$\ are Meir-Keeler, there exist $\delta
_{\varepsilon }>0$\ such that 
\begin{equation}
d(f_{i}(x),f_{i}(y))<\varepsilon \text{,}  \tag{1}
\end{equation}%
for all $i\in I$ and all $x,y\in X^{m}$\ having the property that $d_{\max
}(x,y)<\varepsilon +\delta _{\varepsilon }$. If $%
g=(g^{1},...,g^{m}),h=(h^{1},...,h^{m})\in \mathcal{C}^{m}$ are such that $%
d_{\max }(g,f)<\varepsilon +\delta _{\varepsilon }$, then%
\begin{equation*}
d_{\max }((g^{1}(\alpha (1)),...,g^{m}(\alpha (m))),(h^{1}(\alpha
(1)),...,h^{m}(\alpha (m))))=
\end{equation*}%
\begin{equation*}
=\max \{d(g^{1}(\alpha (1)),h^{1}(\alpha (1))),...,d(g^{m}(\alpha
(m)),h^{m}(\alpha (m)))\}\leq 
\end{equation*}%
\begin{equation*}
\leq \max \{d_{u}(g^{1},h^{1}),...,d_{u}(g^{m},h^{m})=d_{\max
}(g,f)<\varepsilon +\delta _{\varepsilon }
\end{equation*}%
for every $\alpha \in \Omega $. Then, taking into account $(1)$, we get 
\begin{equation*}
d(f_{\alpha ^{1}}(g^{1}(\alpha (1)),...,g^{m}(\alpha (m))),f_{\alpha
^{1}}(h^{1}(\alpha (1)),...,h^{m}(\alpha (m))))=
\end{equation*}%
\begin{equation}
=d(H_{\mathcal{F}}(g)(\alpha ),H_{\mathcal{F}}(h)(\alpha ))<\varepsilon 
\text{,}  \tag{2}
\end{equation}%
for every $\alpha =\alpha ^{1}\alpha ^{2}...\alpha ^{k}\alpha ^{k+1}...\in
\Omega $ and every $g=(g^{1},...,g^{m}),h=(h^{1},...,h^{m})\in \mathcal{C}%
^{m}$ such that $d_{\max }(g,f)<\varepsilon +\delta _{\varepsilon }$. As the
metric space $(\Omega ,d)$\textit{\ }is compact (see Remark 1.3 and take
into account that $I$\ is finite), there exists $\alpha _{0}\in \Omega $
such that 
\begin{equation*}
d_{u}(H_{\mathcal{F}}(g),H_{\mathcal{F}}(h))=\underset{\alpha \in \Omega }{%
\sup }\text{ }d(H_{\mathcal{F}}(g)(\alpha ),H_{\mathcal{F}}(h)(\alpha
))=d(H_{\mathcal{F}}(g)(\alpha _{0}),H_{\mathcal{F}}(h)(\alpha _{0}))\text{.}
\end{equation*}%
In view of $(2)$ we conclude that $d_{u}(H_{\mathcal{F}}(g),H_{\mathcal{F}%
}(h))<\varepsilon $ for every $g,h\in \mathcal{C}^{m}$ such that $d_{\max
}(g,f)<\varepsilon +\delta _{\varepsilon }$, i.e. $H_{\mathcal{F}}$ is
Meir-Keeler. $\square $

\bigskip

\textbf{Proposition 3.6.} \textit{For every comparison function }$\varphi $ 
\textit{and every generalized possibly infinite iterated function system} $%
\mathcal{F}=((X,d),(f_{i})_{i\in I})$ \textit{of order }$m$\textit{\ such
that all the functions }$f_{i}$\textit{\ are }$\varphi $\textit{%
-contractions, the operator }$H_{\mathcal{F}}$ \textit{is Meir-Keeler.}

\textit{Proof}. Let us suppose that the family of functions $(f_{i})_{i\in
I} $\ is not uniformly Meir-Keeler. Then there exists $\varepsilon _{0}>0$
with the property that for every $\delta ,\lambda >0$ there exist $x_{\delta
,\lambda },y_{\delta ,\lambda }\in X^{m}$ and $i_{0}\in I$ such that 
\begin{equation}
d_{\max }(x_{\delta ,\lambda },y_{\delta ,\lambda })<\varepsilon _{0}+\delta
\tag{1}
\end{equation}%
and%
\begin{equation}
d(f_{i_{0}}(x_{\delta ,\lambda }),f_{i_{0}}(y_{\delta ,\lambda }))\geq
\varepsilon _{0}-\lambda \text{.}  \tag{2}
\end{equation}%
Consequently we get $\varepsilon _{0}-\lambda \overset{(2)}{\leq }%
d(f_{i_{0}}(x_{\delta ,\lambda }),f_{i_{0}}(y_{\delta ,\lambda }))\overset{%
f_{i_{0}}\text{ is }\varphi \text{-contraction}}{\leq }$\linebreak $\varphi
(d_{\max }(x_{\delta ,\lambda },y_{\delta ,\lambda })\overset{\text{(1) and
i) from Definition 1.7}}{\leq }\varphi (\varepsilon _{0}+\delta )$, so%
\begin{equation}
\varepsilon _{0}-\lambda \leq \varphi (\varepsilon _{0}+\delta )\text{,} 
\tag{3}
\end{equation}%
for every $\delta ,\lambda >0$. Based on the right continuity of $\varphi $
(see ii) from Definition 1.7), by passing to limit in $(1)$ as $\delta
,\lambda \rightarrow 0$, we get the contradiction $\varepsilon _{0}\leq
\varphi (\varepsilon _{0})\overset{\text{iii) from Definition 1.7}}{<}%
\varepsilon _{0}$.

Hence the family of functions $(f_{i})_{i\in I}$\ is uniformly Meir-Keeler
and Proposition 3.4 assures us that $H_{\mathcal{F}}$ is Meir-Keeler. $%
\square $

\bigskip

\textbf{Proposition 3.7.} \textit{For every generalized possibly infinite
iterated function system} $\mathcal{F}=((X,d),(f_{i})_{i\in I})$ \textit{of
order }$m$\textit{\ such that }$\underset{i\in I}{\sup }$ lip$(f_{i})<1$%
\textit{, the operator }$H_{\mathcal{F}}$ \textit{is Meir-Keeler.}

\textit{Proof}. Indeed, we note that all the functions $f_{i}$\ are $\varphi
_{0}$-contractions, where the comparison function $\varphi _{0}:[0,\infty
)\rightarrow \lbrack 0,\infty )$ is given by $\varphi _{0}(t)=\underset{i\in
I}{(\sup }$ lip$(f_{i}))t$ for every $t\in \lbrack 0,\infty )$. Therefore,
taking into account Proposition 3.6, we conclude that $H_{\mathcal{F}}$ is
Meir-Keeler. $\square $

\bigskip

\textbf{Proposition 3.8.} \textit{For every generalized possibly infinite
iterated function system} $\mathcal{F}=((X,d),(f_{i})_{i\in I})$ \textit{of
order }$m$\textit{\ such that }$I$ \textit{is finite and all the functions} $%
f_{i}$\textit{\ are contractive, the operator }$H_{\mathcal{F}}$ \textit{is
contractive.}

\textit{Proof}. For every $g=(g^{1},...,g^{m}),h=(h^{1},...,h^{m})\in 
\mathcal{C}^{m}$, $g\neq h$, as the metric space $(\Omega ,d)$\textit{\ }is
compact (see Remark 1.3 and take into account that $I$\ is finite), there
exists $\beta \in \Omega $ such that 
\begin{equation*}
d_{u}(H_{\mathcal{F}}(g),H_{\mathcal{F}}(h))=\underset{\alpha \in \Omega }{%
\sup }\text{ }d(H_{\mathcal{F}}(g)(\alpha ),H_{\mathcal{F}}(h)(\alpha
))=d(H_{\mathcal{F}}(g)(\beta ),H_{\mathcal{F}}(h)(\beta ))\text{.}
\end{equation*}%
Then%
\begin{equation*}
d_{u}(H_{\mathcal{F}}(g),H_{\mathcal{F}}(h))=
\end{equation*}%
\begin{equation*}
=d(f_{\beta ^{1}}(g^{1}(\beta (1)),...,g^{m}(\beta (m))),f_{\beta
^{1}}(h^{1}(\beta (1)),...,h^{m}(\beta (m))))\overset{f_{\beta ^{1}}\text{
contractive}}{<}
\end{equation*}%
\begin{equation*}
<d_{\max }((g^{1}(\beta (1)),...,g^{m}(\beta (m)),(h^{1}(\beta
(1)),...,h^{m}(\beta (m))))=
\end{equation*}%
\begin{equation*}
=\max \{d(g^{1}(\beta (1)),h^{1}(\beta (1))),...,d(g^{m}(\beta
(m)),h^{m}(\beta (m)))\}\leq 
\end{equation*}%
\begin{equation*}
\leq \max \{d_{u}(g^{1},h^{1}),...,d_{u}(g^{m},h^{m})=d_{\max }(g,f)\text{,}
\end{equation*}%
for all $g,h\in \mathcal{C}^{m}$, $g\neq h$, where $\beta =\beta ^{1}\beta
^{2}...\beta ^{k}...\in \Omega $, i.e. $H_{\mathcal{F}}$ is contractive. $%
\square $

\bigskip

\textbf{4. The main results}

\bigskip

\textbf{Theorem 4.1.} \textit{For a generalized possibly infinite iterated
function system} $\mathcal{F}=((X,d),(f_{i})_{i\in I})$ \textit{there exists
a unique} $\pi _{0}\in \mathcal{C}$ \textit{such that:}

\textit{a)} 
\begin{equation*}
H_{\mathcal{F}}(\pi _{0},...,\pi _{0})=\pi _{0}\text{;}
\end{equation*}

\textit{b)}%
\begin{equation*}
\underset{n\rightarrow \infty }{\lim }H_{\mathcal{F}}^{[n]}(f,...,f)=\pi _{0}%
\text{,}
\end{equation*}%
\textit{for every} $f\in \mathcal{C}$\textit{, provided that one of the
following conditions is satisfied:}

\textit{i)} \textit{there exists a comparison function }$\varphi $ \textit{%
such that all the functions }$f_{i}$\textit{\ are }$\varphi $\textit{%
-contractions (in particular this happens if }$\underset{i\in I}{\sup }$ lip$%
(f_{i})<1$\textit{);}

\textit{ii)} \textit{the family of functions }$(f_{i})_{i\in I}$\textit{\ is
uniformly Meir-Keeler;}

\textit{iii)} $I$ \textit{is finite and all the functions} $f_{i}$\textit{\
are Meir-Keeler.}

\textit{Proof}. Just use Proposition 3.4, Proposition 3.5, Proposition 3.6
and Theorem 1.11. $\square $

\bigskip

\textbf{Proposition 4.2}. \textit{In the framework of the previous theorem,
we have} $\overline{\pi _{0}(\Omega )}=A_{\mathcal{F}}$.

\textit{Proof}. In view of Theorem 4.1, we have 
\begin{equation*}
\overline{\pi _{0}(\Omega )}\overset{\text{Theorem 4.1, a)}}{=}\overline{H_{%
\mathcal{F}}(\pi _{0},...,\pi _{0})(\Omega )}=\overline{\underset{\alpha \in
\Omega }{\cup }\{H_{\mathcal{F}}(\pi _{0},...,\pi _{0})(\alpha )\}}
\end{equation*}%
\begin{equation*}
=\overline{\underset{\alpha =\alpha ^{1}\alpha ^{2}...\alpha ^{k}\alpha
^{k+1}...\in \Omega }{\cup }\{f_{\alpha ^{1}}(\pi _{0}(\alpha (1)),...,\pi
_{0}(\alpha (m)))\}}=
\end{equation*}%
\begin{equation*}
=\overline{\underset{i\in \Omega }{\cup }f_{i}(\pi _{0}(\Omega ),...,\pi
_{0}(\Omega ))}\overset{(\ast )}{=}\overline{\underset{i\in \Omega }{\cup }%
f_{i}(\overline{\pi _{0}(\Omega )},...,\overline{\pi _{0}(\Omega )})}=
\end{equation*}%
\begin{equation*}
=\mathcal{F}_{\mathcal{F}}(\overline{\pi _{0}(\Omega )},...,\overline{\pi
_{0}(\Omega )})\text{,}
\end{equation*}%
i.e.%
\begin{equation}
\mathcal{F}_{\mathcal{F}}(\overline{\pi _{0}(\Omega )},...,\overline{\pi
_{0}(\Omega )})=\overline{\pi _{0}(\Omega )}\text{,}  \tag{1}
\end{equation}%
where the equality $(\ast )$ is justified in the following way: we have 
\begin{equation*}
f_{i}(\overline{\pi _{0}(\Omega )},...,\overline{\pi _{0}(\Omega )})\overset{%
f_{i}\text{ continuous}}{\subseteq }\overline{f_{i}(\pi _{0}(\Omega
),...,\pi _{0}(\Omega ))}\text{,}
\end{equation*}%
so 
\begin{equation*}
\underset{i\in \Omega }{\cup }f_{i}(\overline{\pi _{0}(\Omega )},...,%
\overline{\pi _{0}(\Omega )})\subseteq \underset{i\in \Omega }{\cup }%
\overline{f_{i}(\pi _{0}(\Omega ),...,\pi _{0}(\Omega ))}\subseteq \overline{%
\underset{i\in \Omega }{\cup }f_{i}(\pi _{0}(\Omega ),...,\pi _{0}(\Omega ))}
\end{equation*}%
and therefore 
\begin{equation*}
\overline{\underset{i\in \Omega }{\cup }f_{i}(\overline{\pi _{0}(\Omega )}%
,...,\overline{\pi _{0}(\Omega )})}\subseteq \overline{\underset{i\in \Omega 
}{\cup }f_{i}(\pi _{0}(\Omega ),...,\pi _{0}(\Omega ))}\subseteq \overline{%
\underset{i\in \Omega }{\cup }f_{i}(\overline{\pi _{0}(\Omega )},...,%
\overline{\pi _{0}(\Omega )})}\text{.}
\end{equation*}

From Remark 1.14, ii) and Theorem 1.15, as $\overline{\pi _{0}(\Omega )}\in 
\mathcal{B}(X)$, using relation $(1)$, we conclude that $\overline{\pi
_{0}(\Omega )}=A_{\mathcal{F}}$. $\square $

\bigskip

\textbf{Proposition 4.3}. \textit{In the framework of the previous theorem,} 
$\pi _{0}$ \textit{is the canonical projection associated to} $\mathcal{F}$.

\textit{Proof}. For $i\in I$, let us consider the function $F_{i}:\Omega
^{m}\rightarrow \Omega $ described in the following way:%
\begin{equation*}
F_{i}(\beta _{1},...,\beta _{m})=\alpha _{1}\alpha _{2}...\alpha _{p}...%
\text{,}
\end{equation*}%
where if, $\beta _{k}=\beta _{1}^{k}\beta _{2}^{k}...\beta _{p}^{k}...\in
\Omega $, with $\beta _{p}^{k}\in \Omega _{p}$, $k\in \{1,...,m\}$, we have $%
\alpha _{1}=i\in \Omega _{1}$, $\alpha _{2}=\beta _{1}^{1}\beta
_{1}^{2}...\beta _{1}^{m}\in \Omega _{2}$, ..., $\alpha _{p+1}=\beta
_{p}^{1}\beta _{p}^{2}...\beta _{p}^{m}\in \Omega _{p+1}$, ... .

\textbf{Claim}. 
\begin{equation}
\pi _{0}(F_{\alpha ^{1}\alpha ^{2}...\alpha ^{k}}(\Lambda _{1},...,\Lambda
_{m^{k}}))=f_{\alpha ^{1}\alpha ^{2}...\alpha ^{k}}(\pi _{0}(\Lambda
_{1}),...,\pi _{0}(\Lambda _{m^{k}}))\text{,}  \tag{1}
\end{equation}%
for all $k\in \mathbb{N}^{\ast }$, $\alpha ^{1}\in I,\alpha ^{2}\in \Omega
_{2},...,\alpha ^{k}\in \Omega _{k}$ and $\Lambda _{1},...,\Lambda
_{m^{k}}\subseteq \Omega $.

\textit{Justification}. We are going to use the mathematical induction
method.

Let us start by noting that%
\begin{equation*}
\pi _{0}(\alpha )\overset{\text{Theorem 4.1, a)}}{=}H_{\mathcal{F}}(\pi
_{0},...,\pi _{0})(\alpha )=f_{\alpha ^{1}}(\pi _{0}(\alpha (1)),...,\pi
_{0}(\alpha (m)))\text{,}
\end{equation*}%
for every $\alpha =\alpha ^{1}\alpha ^{2}...\alpha ^{k}...\in \Omega $, so%
\begin{equation}
\pi _{0}(F_{\alpha ^{1}}(\Lambda _{1},...,\Lambda _{m}))=f_{\alpha ^{1}}(\pi
_{0}(\Lambda _{1}),...,\pi _{0}(\Lambda _{m}))\text{,}  \tag{2}
\end{equation}%
for all $\alpha ^{1}\in I$ and $\Lambda _{1},...,\Lambda _{m}\subseteq
\Omega $, i.e. $(1)$ is valid for $k=1$.

Now we suppose that $(1)$ is valid for $k$ and we prove that is valid also
for $k+1$. Indeed, we have%
\begin{equation*}
\pi _{0}(F_{\alpha ^{1}\alpha ^{2}...\alpha ^{k}\alpha ^{k+1}}(\Lambda
_{1},...,\Lambda _{m^{k+1}}))=
\end{equation*}%
\begin{equation*}
=\pi _{0}(F_{\alpha ^{1}}(F_{\alpha (1)}(\Lambda _{1},...,\Lambda
_{m^{k}}),...,F_{\alpha (m)}(\Lambda _{m^{k+1}-m^{k}+1},...,\Lambda
_{m^{k+1}})))\overset{(2)}{=}
\end{equation*}%
\begin{equation*}
=f_{\alpha ^{1}}(\pi _{0}(F_{\alpha (1)}(\Lambda _{1},...,\Lambda
_{m^{k}})),...,\pi _{0}(F_{\alpha (m)}(\Lambda
_{m^{k+1}-m^{k}+1},...,\Lambda _{m^{k+1}}))))\overset{\text{Claim for }k}{=}
\end{equation*}%
\begin{equation*}
=f_{\alpha ^{1}}(f_{\alpha (1)}(\pi _{0}(\Lambda _{1}),...,\pi _{0}(\Lambda
_{m^{k}})),...,f_{\alpha (m)}(\pi _{0}(\Lambda _{m^{k+1}-m^{k}+1}),...,\pi
_{0}(\Lambda _{m^{k+1}})))=
\end{equation*}%
\begin{equation*}
=f_{\alpha ^{1}\alpha ^{2}...\alpha ^{k}\alpha ^{k+1}}(\pi _{0}(\Lambda
_{1}),...,\pi _{0}(\Lambda _{m^{k+1}}))\text{,}
\end{equation*}%
where $\alpha =\alpha ^{1}\alpha ^{2}...\alpha ^{k}...$, and the
justification of the claim is done.

For $\alpha =\alpha ^{1}\alpha ^{2}...\alpha ^{k}...$ arbitrarily chosen in $%
\Omega $, where $\alpha ^{k}\in I$, we have 
\begin{equation*}
\pi _{0}(\alpha )\in \pi _{0}(\overline{\underset{k\in \mathbb{N}^{\ast }}{%
\cap }F_{\alpha ^{1}\alpha ^{2}...\alpha ^{k}}(\Omega ,...,\Omega )}%
)\subseteq 
\end{equation*}%
\begin{equation*}
\subseteq \overline{\pi _{0}(\underset{k\in \mathbb{N}^{\ast }}{\cap }%
F_{\alpha ^{1}\alpha ^{2}...\alpha ^{k}}(\Omega ,...,\Omega ))}\subseteq 
\overline{\underset{k\in \mathbb{N}^{\ast }}{\cap }\pi _{0}(F_{\alpha
^{1}\alpha ^{2}...\alpha ^{k}}(\Omega ,...,\Omega ))}\overset{\text{Claim}}{=%
}
\end{equation*}%
\begin{equation*}
=\overline{\underset{k\in \mathbb{N}^{\ast }}{\cap }f_{\alpha ^{1}\alpha
^{2}...\alpha ^{k}}(\pi _{0}(\Omega ),...,\pi _{0}(\Omega ))}\subseteq 
\overline{\underset{k\in \mathbb{N}^{\ast }}{\cap }f_{\alpha ^{1}\alpha
^{2}...\alpha ^{k}}(\overline{\pi _{0}(\Omega )},...,\overline{\pi
_{0}(\Omega )})}
\end{equation*}%
\begin{equation*}
\overset{\text{Proposition 4.2}}{=}\overline{\underset{k\in \mathbb{N}^{\ast
}}{\cap }f_{\alpha ^{1}\alpha ^{2}...\alpha ^{k}}(A_{\mathcal{F}},...,A_{%
\mathcal{F}})}\overset{\text{Remark 1.14}}{=}\{\pi (\alpha )\}\text{,}
\end{equation*}%
so $\pi _{0}=\pi $, i.e. $\pi _{0}$ is the canonical projection associated
to $\mathcal{F}$. $\square $

\bigskip

\textbf{5. An example}

\bigskip

Following [2], we consider the generalized possibly infinite iterated
function system $\mathcal{F}=((X,d),(f_{n})_{n\in \mathbb{N}})$ of order $2$%
, where $X=[0,1]$, $d$ is the euclidean metric and the functions $%
f_{n}:[0,1]\times \lbrack 0,1]\rightarrow \lbrack 0,1]$ are given by 
\begin{equation*}
f_{n}(x,y)=\frac{1}{2^{n+2}}(x+y)+\frac{1}{2^{n+1}}\text{,}
\end{equation*}%
for every $n\in \mathbb{N}$, $x,y\in \lbrack 0,1]$.

Note that:

i) the family $(f_{n})_{n\in \mathbb{N}}$ is uniformly Meir-Keeler since for
every $\varepsilon >0$ there exist $\lambda _{\varepsilon }=\frac{%
\varepsilon }{9}>0$ and $\delta _{\varepsilon }=\frac{\varepsilon }{3}>0$
such that $\left\vert f_{n}(u)-f_{n}(v)\right\vert <\varepsilon -\lambda
_{\varepsilon }$ for all $u,v\in \lbrack 0,1]\times \lbrack 0,1]$ with the
property that $d_{\max }(u,v)<\varepsilon +\delta _{\varepsilon }$;

ii) the operator $H_{\mathcal{F}}:\mathcal{C}^{2}\rightarrow \mathcal{C}$
acts in the following way:%
\begin{equation*}
H_{\mathcal{F}}(g_{1},g_{2})(\alpha )=\frac{1}{2^{\alpha ^{1}+2}}%
(g_{1}(\alpha _{2}^{1}\alpha _{3}^{1}...\alpha _{k}^{1}...)+g_{2}(\alpha
_{2}^{2}\alpha _{3}^{2}...\alpha _{k}^{2}...))+\frac{1}{2^{\alpha ^{1}+1}}%
\text{,}
\end{equation*}%
for every $g_{1},g_{2}\in \mathcal{C}$ and every $\alpha =\alpha ^{1}\alpha
^{2}...\alpha ^{k}...\in \Omega $, where $\alpha ^{k}=\alpha _{k}^{1}\alpha
_{k}^{2}\in \Omega _{k}$ with $\alpha _{k}^{1},\alpha _{k}^{2}\in \Omega
_{k-1}$;

iii) $A_{\mathcal{F}}=[0,1]$.

\bigskip

\textbf{References}

\bigskip

[1] D. Dumitru, Generalized iterated function systems containing Meir-Keeler
functions, An. Univ. Bucur., Mat., \textbf{58 }(2009), 109-121.

[2] D. Dumitru, Contraction-type functions and some applications to
\linebreak GIIFS, to appear in An. Univ. Spiru Haret, Ser. Mat.-Inform.

[3] J.E. Hutchinson, Fractals and self similarity, Indiana Univ. Math. J., 
\textbf{30 }(1981), 713-747.

[4] P. Jaros, \L . Ma\'{s}lanka and F. Strobin, Algorithms generating images
of attractors of generalized iterated function systems, Numer. Algorithms, 
\textbf{73} (2016), 477-499.

[5] R. Miculescu, Generalized iterated function systems with place dependent
probabilities, Acta Appl. Math., \textbf{130} (2014), 135-150.

[6] A. Mihail and R. Miculescu, Applications of Fixed Point Theorems in the
Theory of Generalized IFS, Fixed Point Theory Appl. Volume 2008, Article ID
312876, 11 pages doi: 10.1155/2008/312876.

[7] A. Mihail and R. Miculescu, A generalization of the Hutchinson measure,
Mediterr. J. Math., \textbf{6} (2009), 203--213.

[8] A. Mihail and R. Miculescu, Generalized IFSs on Noncompact Spaces ,
Fixed Point Theory Appl. Volume 2010, Article ID 584215, 11 pages doi:
10.1155/2010/584215.

[9] A. Mihail, The shift space for recurrent iterated function systems, Rev.
Roum. Math. Pures Appl., \textbf{53} (2008), 339-355.

[10] A. Mihail, The canonical projection between the shift space of an IIFS
and its attractor as a fixed point, Fixed Point Theory Appl., 2015, Paper
No. 75, 15 p.

[11] E. Oliveira and F. Strobin, Fuzzy attractors appearing from GIFZS,
Fuzzy Set Syst., in print, doi:10.1016/j.fss.2017.05.003.

[12] N.A. Secelean, Invariant measure associated with a generalized
countable iterated function system, Mediterr. J. Math., \textbf{11} (2014),
361-372.

[13] N.A. Secelean, Generalized iterated function systems on the space $%
l^{\infty }(X)$, J. Math. Anal. Appl., \textbf{410} (2014), 847-858.

[14] F. Strobin, Attractors of generalized IFSs that are not attractors of
IFSs, J. Math. Anal. Appl., \textbf{422 }(2015), 99-108.

[15] F. Strobin and J. Swaczyna, On a certain generalisation of the iterated
function system, Bull. Aust. Math. Soc., \textbf{87} (2013), 37-54.

[16] F. Strobin and J. Swaczyna, A code space for a generalized IFS, Fixed
Point Theory, \textbf{17} (2016), 477-493.

\bigskip

{\small Radu MICULESCU}

{\small Faculty of Mathematics and Computer Science}

{\small University of Bucharest, Romania}

{\small Academiei Street 14, 010014, Bucharest, Romania}

{\small E-mail: miculesc@yahoo.com}

\bigskip

{\small Silviu-Aurelian URZICEANU}

{\small Faculty of Mathematics and Computer Science}

{\small University of Pite\c{s}ti, Romania}

{\small T\^{a}rgul din Vale 1, 110040, Pite\c{s}ti, Arge\c{s}, Romania}

{\small E-mail: fmi\_silviu@yahoo.com}

\newpage

\end{document}